\documentclass[12pt]{article}

\RequirePackage{amssymb}
\RequirePackage{amsfonts}
\RequirePackage{amsmath}
\RequirePackage{graphicx,epstopdf} 
\RequirePackage{array}
\RequirePackage[mathcal]{euscript}
\RequirePackage{amsthm}
\RequirePackage{latexsym,amscd}
\RequirePackage{url}
\RequirePackage{dsfont}
\RequirePackage{lmodern} 
\RequirePackage{caption}
\RequirePackage{forloop}
\RequirePackage{hyperref}

\hypersetup{bookmarksdepth = 0}  

\newcommand{\E}{\mathbb{E}}

\renewcommand{\P}{\mathbb{P}}

\newcommand{\R}{\mathbb{R}}

\newcommand{\Z}{\mathbb{Z}}

\newcommand{\cN}{\mathcal{N}}

\newcommand{\cS}{\mathcal{S}}

\renewcommand{\th}{\theta}

\renewcommand{\l}{\lambda}

\newcommand{\s}{\sigma}

\newcommand{\ds}{\: \mathrm{d}s}

\newcommand{\dx}{\: \mathrm{d}x}

\newcommand{\rmd}{\mathrm{d}}

\newcommand{\pinf}{+ \infty}
\newcommand{\eqlaw}{\overset{(d)}{=}}

\newcommand{\un}{\mathds{1}}
\newcommand{\unn}[1]{\mathds{1}_{\left \{ #1 \right \} }}

\RequirePackage[english]{babel}
\RequirePackage{vmargin}
\RequirePackage{enumitem}

\hypersetup{colorlinks,citecolor=green,filecolor=purple,linkcolor=red,urlcolor=blue} 

\setmarginsrb{2cm}{2cm}{2cm}{3cm}{0cm}{0cm}{0cm}{1cm}

\newtheorem{thm}{Theorem}[section]
\newtheorem{theorem}[thm]{Theorem}
\newtheorem{prop}[thm]{Proposition}

\newtheorem{cor}[thm]{Corollary}

\theoremstyle{definition}

\theoremstyle{remark}
\newtheorem{remark}[thm]{Remark}

\makeatletter

\newcounter{i}
\newcounter{aut}
\newcounter{loc}

\renewcommand{\author}[2][1]{\addtocounter{aut}{1} \global \@namedef{type@\arabic{aut}}{#1} \@namedef{author@\arabic{aut}}{#2}}
\newcommand{\location}[2][1]{\addtocounter{loc}{1} \global \@namedef{location@\arabic{loc}}{#2}}
\def \gettype#1{\@nameuse{type@\arabic{#1}}}
\def \getauthor#1{\@nameuse{author@\arabic#1}}
\def \getlocation#1{\@nameuse{location@\arabic#1}}
\def \email{\@namedef{@email}}

\def \maketitle
   {%
   \renewcommand{\thefootnote}{\fnsymbol{footnote}}%
   \begin{center}%
    \LARGE \@title%
    \vskip 0.6cm%
    \forloop{i}{1}{\value{i} < \value{aut}}%
       {%
       \large \getauthor{i}\footnotemark[\gettype{i}] \hspace{1cm} %
       }%
      \large \@nameuse{author@\arabic{aut}}%
       \footnotemark[\@nameuse{type@\arabic{aut}}]%
     \forloop{i}{1}{\value{i} < \value{loc} \or \value{i} = \value{loc}}%
       {%
       \footnotetext[\arabic{i}]{\getlocation{i}}
       }%
     \footnotetext%
       {%
       E-mail address: \@nameuse{@email}%
       }%
    \end{center}%
	 \renewcommand{\thefootnote}{\arabic{footnote}} \setcounter{footnote}{0}%
  }%

\makeatother

\newcommand{\keywords}{\vspace*{0.3cm} \noindent \textbf{Keywords}: }
\newcommand{\MSC}{\vspace*{0.3cm} \noindent \textbf{MSC 2010}: }

\title{Local times in a Brownian excursion}
\author[1]{Krishna B. \textsc{Athreya}}
\author[2]{Raoul \textsc{Normand}}
\author[1]{Vivekananda \textsc{Roy}}
\author[3]{Sheng-Jhih \textsc{Wu}}
\location[1]{Iowa State University, Ames, Iowa 50011, USA}
\location[2]{Institute of Mathematics, Academia Sinica, Taipei 10617, Taiwan}
\location[3]{Center for Advanced Statistics and Econometrics Research, School of Mathematical Sciences, Soochow University, Suzhou 215006, China}
\email{kbathreya@gmail.com; rnormand@math.sinica.edu.tw; vroy@iastate.edu; \\ szwu@suda.edu.cn}

\begin{document}

\maketitle

\begin{abstract}
  Let $\{B(t), t \geq 0\}$ be a standard Brownian motion in $\R$. Let
  $T$ be the first return time to 0 after hitting 1, and $\{L(T,x), x
  \in \R\}$ be the local time process at time $T$ and level $x$. The
  distribution of $L(T,x)$ for each $x \in \R$ is determined. This is
  applied to the estimation of a $L^1$ integral on $\R$.
\end{abstract}

\keywords Brownian motion, local time, Ray-Knight theorem, regenerative process.

\MSC Primary: 60J65; Secondary: 60F05

\section{Introduction}

Let $B \equiv \{B(t), t \geq 0\}$ be the standard Brownian motion (SBM) on
$\R$, and $\{L(t,x), t \geq 0, x \in \R \}$ be its {\it local time
  process}. That is, with probability 1, for any Borel set $A \subset
\R$ and $t \geq 0$,
\begin{equation}
  \label{eq:defloct}
  \int_0^t \un_{A}(B(s)) \ds = \int_A L(t,x) \dx.
\end{equation}
Such a process exists (see e.g. Chapter VI in \cite{RY}). Let
\begin{equation}
  \label{eq:defregen}
  T = \inf \{ t > 0, \; B(t) = 0 \; \mbox{and for some} \; s \in (0,t), \;B(s) = 1 \}
\end{equation}
be the first hitting time of 0 after hitting 1. It may be noted that
by the strong Markov property $\{B(T+u), u \geq 0\}$ is independent of
$\{B(u), 0 \le u < T \}$ and has the same distribution as SBM. That
is, the process {\it regenerates } at time $T$. This and similar ideas have been used in
\cite{AthreyaDKR} to establish limit theorems for Brownian motion
and stable processes.  The current paper is devoted to finding the
marginal distribution of $L(T,x)$ for $x \in \R$. This is given in
Theorem \ref{th:Ltx}, along with some useful corollaries. The proofs
are in Section \ref{sec:proofs}. An application (and the motivation
for this work) is to obtain a point estimate and a confidence interval
based on $\{B(t), t \geq 0\}$ for $\int_{\R} f(x) \dx$, where $f \in
L^1(\R)$. This is done in Section \ref{sec:application}, while Section
\ref{sec:open} gives an extension of the main result and presents an
open problem.

\section{Main results}

In the following, we define, for $x \in \R$,
\[
\tau_x = \inf \{ t \geq 0, B(t) = x \}
\]
to be the hitting time of $x$ for SBM. Below $X \eqlaw Y$ means $X$ and
$Y$ are random variables with the same distribution. We establish the following.

\begin{theorem} \label{th:Ltx}

  Let $Y_1$ and $Y_2$ be two independent chi-squared random variables with
  two degrees of freedom, i.e. $Y_1 \eqlaw Y_2 \eqlaw Z_1^2 + Z_2^2$
  for independent standard normal $Z_1$ and $Z_2$. Let $T$ be as in \eqref{eq:defregen}. Assume that
  $Y_1$ and $Y_2$ are independent of $B$. Then, the following hold:
\begin{enumerate}
\item[(i)] for $x > 1$, $L(T,x) \eqlaw \unn{\tau_x \leq T} x Y_1$,
\item[(ii)] for $x \in [0,1]$, $L(T,x) \eqlaw x Y_1 + (1-x)Y_2$,
\item[(iii)] for $x < 0$, $L(T,x) \eqlaw \unn{\tau_x \leq \tau_1} (1-x) Y_1$,
\end{enumerate}
where $\{L(t, x), t \ge 0, x \in \mathbb{R}\}$ is as in \eqref{eq:defloct}.
\end{theorem}

\begin{cor} \label{cor:expmoment} For all $x \in \R$, there exists
  $\theta(x) >0$, such that $\E[\exp \{\th L(T,x)\}] < \infty$ for all
  $\theta < \theta(x)$. Hence, $L(T,x)$ has all moments finite.
\end{cor}

\begin{cor} \label{cor:moments}
The first and second moments of  $L(T,x)$ are given below.
\begin{enumerate}
  \item [(i)] For all $x \in \R$, $\E(L(T,x)) = 2$, and
\item [(ii)] \[
\E \left ( L(T,x)^2 \right ) =
\begin{cases}
8x & \text{for $x > 1$}, \\
8(x^2 - x + 1) & \text{for $x \in [0,1]$}, \\
8(1-x) & \text{for $x < 0$}. \\
\end{cases}
\]
  \end{enumerate}
\end{cor}

\begin{cor} \label{cor:Ef}
Let $f : \R \to \R$ be Lebesgue integrable, and assume that $\int_{\R} \sqrt{|x|} |f(x)| \dx < \infty$. Then
\[
\E \left ( \int_0^T \left | f(B(s)) \right | \ds \right )^2 < \infty,
\]
where $T$ is as in \eqref{eq:defregen} and $\{B(s), s \ge 0\}$ is SBM.
\end{cor}

\section{Proofs} \label{sec:proofs}

\subsection{Proof of Theorem \ref{th:Ltx}}

In the following, we will also consider Brownian motions $B^a \equiv
\{B^a(t) \equiv B(t)+a, t \ge 0\}$ where $\{B(t), t \ge 0\}$ is
SBM. We define $\tau_x^a$ as the hitting time of $x$ by $B^a$. Let
$\{L^a(t,x), t \ge 0, x \in \mathbb{R}\}$ be its local time process.

 Let $x > 1$. On the event $\{\tau_x >
  T\}$, $L(T,x) = 0$. By the strong Markov
  property of $B$, on the event $\{\tau_x \le T\}$,
\[
\{B(\tau_x + t) , t \in [0,T-\tau_x]\} \eqlaw \{B^x(t), t \in [0,\tau^x_{0})]\},
\]
and further, the left side is independent of $\{B(u), 0 \le u \le
\tau_x\}$. Since $B$ does not accumulate any local time at $x$ before
$\tau_x$, i.e. $L(\tau_x,x) = 0$, we have
\[
L(T,x) \eqlaw L^x(\tau^x_0, x).
\]
By translation invariance of SBM, the right side above has the same
distribution as $L(\tau_{-x},0)$. This, in turn by symmetry of SBM
about 0, has the same distribution as $L(\tau_x,0)$. Next, the scaling
property of local times (see exercise 2.11 Chapter VI and Theorem
2.2 Chapter XI in \cite{RY}) ensures that
\[
L(\tau_x,0) \eqlaw x L(\tau_1,0).
\]
Further, the first Ray-Knight theorem (p. 454 in \cite{RY}) asserts
that the process $\{L(\tau_1, 1-a), 0 \le a \le 1 \}$ is a squared Bessel
2 process, i.e. has the distribution same as that of
\begin{equation}
  \label{eq:bess}
\{H(a), 0 \le a \le 1\} \equiv  \left \{ B_1(a)^2 + B_2(a)^2, a \in [0,1] \right \},
\end{equation}
where $B_1$ and $B_2$ are independent standard Brownian motions. In particular,
\[
L(\tau_1, 0) \eqlaw B_1(1)^2 + B_2(1)^2 \eqlaw Y_1.
\]
We can thus conclude that, on $\{\tau_x \le T\}$ ,
\[
L(T,x) \eqlaw x Y_1.
\]
This proves $(i)$.

 Next, let $x \in [0,1]$. The continuity of Brownian trajectories implies that
  $\tau_x \leq \tau_1 \leq T$. Thus the local time at $x$ accumulated on
  $[0,T]$ is the sum of
\begin{itemize}
\item the local time $L(\tau_1,x)$ at $x$ accumulated on $[0,\tau_1]$,
\item the local time accumulated at $x$ on $[\tau_1,T]$.
\end{itemize}
But by the strong Markov property again, the second one is independent
of $L(\tau_1,x)$ and has the distribution of $L^1(\tau^1_0,x)$. By
translation invariance of SBM, this is distributed as
$L(\tau_{-1},x-1)$, which by symmetry of SBM has
same distribution as $L(\tau_1,1-x)$. So,
\begin{equation} \label{eq:LTx}
L(T,x) \eqlaw L_1(\tau_1,x) + L_2(\tau_1,1-x)
\end{equation}
for $L_1, L_2$ two independent copies of the process $L$. By the first
Ray-Knight theorem again, we conclude
\[
L(T,x) \eqlaw (1-x) Y_1 + x Y_2 .
\]
This proves $(ii)$.

Finally, let $x < 0$. On the event $\{\tau_x > \tau_1\}$, we
  have $\tau_x > T$ by continuity of Brownian trajectories, and thus $L(T,x) = 0$. On the other
  hand, on $\{\tau_x \leq \tau_1\}$, we have, as before,
\[
L(T,x) = L(\tau_1,x) \eqlaw L^x(\tau^x_1,x) \eqlaw L(\tau_{1-x},0) \eqlaw (1-x) L(\tau_1,0) \eqlaw (1-x) Y_1,
\]
where the second equality is by the strong Markov property, the third by the
translation invariance of Brownian motion, the fourth by the
scaling property and the fifth by the Ray-Knight theorem. This proves
$(iii)$.

\subsection{Proof of Corollary \ref{cor:expmoment}}
Let $x >1$. From Theorem~\ref{th:Ltx} $(i)$ we know that $L(T,x)
\eqlaw \unn{\tau_x \leq T} x Y_1$. For $\theta >0$ it then follows
that
\[\E[\exp \{\th L(T,x)\}] \le \E[\exp \{\th x Y_1 \}] = \E[\exp
\{\th x (Z_1^2 + Z_2^2) \}],\] where $Z_1, Z_2$ independent $N (0, 1)$
random variables. The proof follows since $\E[\exp \{\th x (Z_1^2 +
Z_2^2) \}] < \infty$ if $\theta < \theta(x) \equiv 1/(2x)$.
Similarly, for $x < 0$, $\E[\exp \{\th L(T,x)\}] < \infty$ if $\theta <
1/[2(1-x)]$ and $\E[\exp \{\th L(T,x)\}] < \infty$ if $\theta <
1/2$ when $x \in [0, 1]$.

\subsection{Proof of Corollary \ref{cor:moments}}

In the following, recall that $\E(Y_1) = 2$ and $\E(Y_1^2) = 8$, since $Y_1 \eqlaw Z_1^2 + Z_2^2$ for
$Z_1, Z_2$ independent $N (0, 1)$ random variables.

Now, for $x > 1$, by Theorem~\ref{th:Ltx} $(i)$, and independence of $\un_{\{\tau_x \le T\}}$ and $Y_1$,
\[
\E(L(T,x)) = \P(\tau_x \le T) x \E(Y_1).
\]
But $\tau_x \le T$ if and only if the Brownian motion hits 1, then
hits $x$ before hitting 0, and thus, by the strong Markov property of
SBM,
\[
\P(\tau_x \le T) = \P(\tau^1_x < \tau^1_0) = \frac1x,
\]
where the last equality is the continuous analog of the usual
gambler's ruin estimate (see Proposition 2.8 Chapter II in \cite{RY}). Thus,
\[
\E(L(T,x)) = \frac1x \times x \times 2 = 2.
\]
Next, for $x \in [0,1]$, by Theorem~\ref{th:Ltx} $(ii)$,
\[
\E(L(T,x)) = x \E(Y_1) + (1-x) E(Y_2) = 2 x + 2(1-x) = 2.
\]
Finally, for $x < 0$, by gambler's ruin estimates again, we have
\[
\P(\tau_x \leq \tau_1) = \frac{1}{1-x}
\]
and thus by Theorem~\ref{th:Ltx} $(iii)$ and independence of
$\un_{\{\tau_x \le \tau_1\}}$ and $Y_1$, we have
\[
\E(L(T,x)) = \frac{1}{1-x} \times (1-x) \times 2 = 2.
\]

\noindent For the second moment, note that, by independence as above
we have by Theorem~\ref{th:Ltx} $(i)$, for $x > 1$,
\[
\E(L(T,x)^2) = \E(\unn{\tau_x \le T}^2) x^2 E(Y_1^2) = \P(\tau_x \le
T) x^2 \E(Y_1^2) = \frac1x \times x^2 \times 8 = 8x.
\]
The same works for $x < 0$. Finally, for $x \in [0,1]$, by independence,
\begin{align*}
\E(L(T,x)^2) & = x^2 E(Y_1^2) + (1-x)^2 E(Y_2^2) + 2 x(1-x) \E(Y_1)\E(Y_2) \\
& = 8 x^2 + 8 (1-x)^2 + 8 x(1-x) \\
& = 8 (x^2 - x + 1)
\end{align*}
and the proof is complete.

\subsection{Proof of Corollary \ref{cor:Ef}}

Since $f$ is integrable, by the definition of $L(T, x)$ as in
\eqref{eq:defloct} (see also Corollary 1.6 Chapter VI in \cite{RY}),
\[
\E \left ( \int_0^T |f(B_s)| \ds \right ) = \int_{\R} |f(x)|
\E(L(T,x)) \dx = 2 \int_{\R} |f(x)|\dx < \infty,
\]
so that $\int_0^T |f(B_s)| \ds < \infty$ a.s. and $V := \int_0^T
f(B_s) \ds$ is a well-defined random variable. Next, by Minkowski
inequality \cite{AthreyaLahiri}, for any $p \in [1,\infty)$,
\[
(\E|V|^p)^{1/p} \leq \int_{\R} |f(x)| (\E L(T,x)^p)^{1/p} \dx.
\]
But, by Corollary \ref{cor:moments} $(ii)$, there is a $C \in
(0,\infty)$ such that for all $x$ in $\R$,
\[
(\E L(T,x)^2)^{1/2} \leq C \sqrt{|x|}.
\]
This yields
\[
(\E |V|^2)^{1/2} \leq C \int_{\R} |f(x)| \sqrt{|x|} \dx,
\]
which is finite by hypothesis.

\begin{remark}
By the same proof, and with an easy extension of Corollary \ref{cor:moments}, it holds that $E(|V|^p) < \infty$ when
\[
\int_{\R} |f(x)| |x|^{\frac{p-1}{p}} \dx < \infty.
\]
In particular, $V$ has moments of all order when
\[
\int_{\R} |f(x)| |x| \dx < \infty.
\]
\end{remark}

\section{An application} \label{sec:application}

Let $f : \R \to \R$ be integrable with respect to the Lebesgue measure
and $\l := \int_{\R} f(x) \dx$. Consider the problem of obtaining a
point estimate and a confidence interval for $\l$ based on a suitable
statistical data.

This problem was solved on $\R^d$ by \cite{AthreyaRoyEI} via simple
symmetric random walks on $\Z$ and an appropriate randomization around
the observed values of the walk. In that paper, the authors introduced
a new Monte Carlo procedure called Regenerative Sequence Monte Carlo
(RSMC), which works for estimating $\int_S f \; \rmd \pi$ for an
integrable function $f$ on a measure space $(S,\cS, \pi)$. While the
classical i.i.d.\ Monte Carlo method and the currently popular Markov Chain
Monte Carlo method require $\pi(S) < \infty$, the RSMC allows $\pi(S)$
to be finite or infinite. The key requirement imposed in
\cite{AthreyaRoyEI} is that the regenerative sequence used be such
that its occupation measure coincides with a constant multiple of the
given measure $\pi$. The regenerative sequence need not even be
Markovian.

Here, we produce a solution similar to that of \cite{AthreyaRoyEI},
but using the regenerative property of standard Brownian motion. This
was exploited by \cite{AthreyaDKR} in proving many
results similar to those of Darling and Kac \cite{DarlingKac}, in
particular giving an alternative proof of a well-known result of
Kallianpur and Robbins \cite{KallianpurRobbins}.

We start with some basic results on standard Brownian motion. For
proofs of the five following results, see \cite{AthreyaDKR}. In the
following, $B$ is still a standard Brownian motion, and $T$ the first
hitting time of 0 after hitting 1 as defined in \eqref{eq:defregen}.

\begin{prop}
It holds that $T < \infty$ a.s. and there is a $0 < C < \infty$ such that
\[
\P(T > y) \sim \frac{C}{\sqrt{y}}, \quad \mbox{as}\;\; y \to \infty.
\]
\end{prop}

\begin{prop}
Let $T_0 = 0$, $T_1 = T$ and for $i \geq 1$,
\[
T_{i+1} = \inf \{ t > T_i, \exists \; s \in (T_i, t)\; \mbox{such
  that}\; \; B(s) = 1\; \mbox{and}\; B(t) = 0\}.
\]
Define $\eta_i = \{ B(t), \; T_{i-1} \leq t \leq T_i, \; T_i - T_{i-1}
\}$ for $i \geq 0$. Then the excursions $\{\eta_i\}_{i \geq 0}$ are i.i.d.
\end{prop}

\begin{prop}
Let
\[
\pi(A) \equiv \E \left ( \int_0^T \un_A(B(s)) \ds \right )
\]
for any Borel set $A \subset \R$. Then, there exists $c \in (0,\infty)$
such that, for any Borel set $A$,
\[
\pi(A) = c \; m(A),
\]
where $m$ is the Lebesgue measure.
\end{prop}

\begin{prop}
Assume that $f \in L^1(\R)$ and define $\l = \int_{\R} f(x) \dx$. Let
\[
\l(t) \equiv \frac{\int_0^t f(B(s)) \ds}{\int_0^t \unn{[0,1]} (B(s))
  \ds}, \quad t > 0.
\]
Then $\l(t) \to \l$ a.s. at $t \to \infty$.
\end{prop}

The next result is proved in \cite{AthreyaDKR} and \cite{KallianpurRobbins}.

\begin{prop}
Let $f \in L^1(\R)$ and
\[
\widetilde{\l}(t) = \frac{1}{\sqrt{t}} \int_0^t f(B(s)) \ds, \quad t > 0.
\]
Then
\[
\widetilde{\l}(t) \overset{d}{\to} \left ( \int_{\R} f(x) \dx \right
) \sqrt{|Z|} \sqrt{\frac{2}{\pi}}
\]
where $Z$ is a $N(0, 1)$ random variable.
\end{prop}

Let $\{N(t), t \geq 0\}$ be the renewal process generated by the
random walk $\{T_i, \; i \geq 0\}$, i.e. $N(t) = k$ for $T_k \leq t <
T_{k+1}$. Then, it follows from the recent work of
\cite{AthreyaRoyEI} that, if $f \in L^1(\R)$ and
\begin{equation} \label{eq:cond}
\E \left ( \int_0^T f(B(s)) \ds \right )^2 < \infty
\end{equation}
then
\[
\frac{\sqrt{N(t)}}{\s} \left ( \frac{1}{N(t)} \int_0^t f(B(s)) \ds - \l \right ) \overset{d}{\to} \cN(0,1)
\]
where
\[
\s^2 =  \E \left ( \int_0^T f(B(s)) \ds \right )^2 - \l^2.
\]
By Corollary \ref{cor:Ef}, a sufficient condition for \eqref{eq:cond}
is that $\int_{\R} |f(x)| \sqrt{|x|} \dx < \infty$. This yields the
following result.

\begin{theorem}
\label{thm:bmregest}
  Let $f : \R \to \R$ be integrable with respect to the Lebesgue
  measure and $\l = \int_{\R} f(x) \dx$. Assume further that
  $\int_{\R} |f(x)| \sqrt{|x|} \dx < \pinf$, and define $N(t)$ and
  $\s^2$ as above, and
  \begin{equation}
    \label{eq:bmregest}
    \l^*(t) \equiv \frac{1}{N(t)} \int_0^t f(B(s)) \ds.
  \end{equation}
Then
\begin{enumerate}
\item[(i)] \hspace{1.7in} $\l^*(t) \to \l$ a.s.,
\item[(ii)]
  \begin{equation*}
    \label{eq:bmregclt}
    \frac{(\l^*(t) - \l)}{\s} \sqrt{N(t)}\overset{d}{\to} \cN(0,1), \quad t \to \infty ,
  \end{equation*}
\item[(iii)]
\[
\frac{(\l^*(t) - \l)}{\s} t^{1/4} \overset{d}{\to} Q, \quad t \to \infty,
\]
where $Q$ is a random variable having the same distribution as of
$B(V)$, where $V$ is a random variable independent of $B$ and has a
distribution of a stable law of order $1/2$, i.e. $\E(e^{-s V}) = e^{-
  \sqrt{s}}$ for $s \geq 0$.
\end{enumerate}
\end{theorem}

The proof is very similar to that of Theorem 1 in \cite{AthreyaRoyEI}
and is omitted.

\begin{remark}
  \label{rem:bmest}
  Now, as in \cite{AthreyaRoyEI} based on Theorem \ref{thm:bmregest}
  $(i)$ and $(ii)$ and the data $\{ B(u), 0 \le u \le t\}$ a
  point estimate for $\lambda$ is provided by $\l^*(t)$ defined in
  \eqref{eq:bmregest} and an asymptotic $(1 - \alpha)$ level $(0 <
  \alpha <1)$ confidence interval for $\lambda$ is $I_t \equiv
  (\l^*(t) - \hat{\sigma} z_\alpha / \sqrt{N(t)}, \l^*(t) +
  \hat{\sigma} z_\alpha / \sqrt{N(t)})$ where $P(|Z| > z_\alpha) =
  \alpha$, $Z \sim N(0, 1)$, $\hat{\sigma}^2 = \sum_{j=1}^{N(t)}
  \xi_j^2 /N(t) - \l^*(t)^2$ and $\xi_j = \int_{T_{j-1}}^{T_j} f(B(u))
  du$, $j \ge 1$.
\end{remark}
\section{A related result and an open problem} \label{sec:open}

From the proof of Theorem \ref{th:Ltx} $(ii)$, more precisely Equation
\eqref{eq:LTx}, we have the following result.

\begin{prop}
Let $H_1 (\cdot)$ and $H_2 (\cdot)$ be two independent Bessel 2
processes as defined in \eqref{eq:bess}. Then
\[
\{L(T,x), x \in [0,1]\} \eqlaw \{H_1(x) + H_2(1-x), x \in [0,1]\}.
\]
\end{prop}

This gives the distribution of the process $\{L(T,x), x \in
[0,1]\}$. A natural open problem is determining the distribution of the whole
process $\{L(T, x), x \in \R \}$.


\begin{thebibliography}{20}

\bibitem{AthreyaDKR}
K.~B. Athreya.
\newblock Darling and {K}ac revisited.
\newblock {\em Sankhy\=a Ser. A}, 48(3):255--266, 1986.

\bibitem{AthreyaLahiri}
K.~B. Athreya and S.~N. Lahiri.
\newblock {\em Measure theory and probability theory}.
\newblock Springer Texts in Statistics. Springer, New York, 2006.

\bibitem{AthreyaRoyEI}
K.~B. Athreya and V.~Roy.
\newblock Estimation of integrals with respect to infinite measures using
  regenerative sequences.
\newblock {\em J. Appl. Probab.}, 52, 2015.
\newblock To appear.

\bibitem{DarlingKac}
D.~A. Darling and M.~Kac.
\newblock On occupation times for {M}arkoff processes.
\newblock {\em Trans. Amer. Math. Soc.}, 84:444--458, 1957.

\bibitem{KallianpurRobbins}
G.~Kallianpur and H.~Robbins.
\newblock Ergodic property of the {B}rownian motion process.
\newblock {\em Proc. Nat. Acad. Sci. U. S. A.}, 39:525--533, 1953.

\bibitem{RY}
D.~Revuz and M.~Yor.
\newblock {\em Continuous martingales and {B}rownian motion}, volume 293 of
  {\em Grundlehren der Mathematischen Wissenschaften [Fundamental Principles of
  Mathematical Sciences]}.
\newblock Springer-Verlag, Berlin, third edition, 1999.

\end{thebibliography}
\end{document}